\documentclass{amsart}
\input epsf

\usepackage{amsfonts,amsthm,amsmath,amssymb,latexsym}
\usepackage{graphicx,color}
\usepackage[all]{xy}

\def\D{{\bf D}}

\def\G{\mathcal{G}}

\begin{document}

\author{Dragomir \v Sari\' c}

\address{Institute for Mathematical Sciences, Stony Brook University,
Stony Brook, NY 11794-3660} \email{saric@math.sunysb.edu}

\title{Bounded earthquakes}


\date{\today}

\begin{abstract}
We give a short proof of the fact that bounded earthquakes of the
unit disk induce quasisymmetric maps of the unit circle. By a
similar method, we show that symmetric maps are induced by bounded
earthquakes with asymptotically trivial measures.
\end{abstract}

\maketitle

\section{Introduction}

An earthquake $E$ of the unit disk $\D$ is a piecewise isometric
(for the hyperbolic metric on $\D$) surjective self-map. The unit
disk $\D$ is partitioned into non-intersecting geodesics and
complementary ideal polygons, called the {\it strata} of $E$; $E$
is an isometry on each stratum. The set of geodesics in the strata
is called the {\it support} of $E$ and it forms a geodesic
lamination of $\D$. By the definition, an earthquake $E$ moves its
strata to the left relative to each other. Thurston \cite{Th}
introduced earthquakes of the unit disk and showed a remarkable
fact that each homeomorphism of the unit circle $S^1$ can be
obtained as an extension to the boundary $\partial\D=S^1$ of an
earthquake map $E:\D\to\D$ \cite[Theorem 3.1]{Th}.

\vskip .2 cm

To an earthquake $E$ corresponds a unique transverse measure to
its support (invariant under isotopies relative the support
geodesic lamination) given by measuring the relative movement to
the left of the strata of $E$ \cite{Th}. Thurston \cite{Th}
defined the norm of an earthquake measure to be the supremum over
all geodesic arcs $I$ of length $1$ in $\D$ of the deposited
transverse measures to $I$. If the norm is finite then the
earthquake measure is called {\it bounded}. The following theorem
is very suggestive from Thurston's paper \cite{Th}:

\vskip .2 cm

\paragraph{\bf Theorem 1} {\it Let $\mu$ be an earthquake measure
which defines an earthquake map $E_{\mu}$ whose extension to $S^1$
is a homeomorphism. Then the following are equivalent:}
\begin{enumerate}

\item $\mu$ is a bounded earthquake measure

\item $E_{\mu}|_{S^1}:S^1\to S^1$ is a quasisymmetric map

\end{enumerate}

\vskip .2 cm

The first proof of the equivalence of the two conditions \cite{Sa}
involved a third condition:

\vskip .2 cm

\hskip .2 cm (3) $t\mapsto E_{t\mu}|_{S^1}(x)$, for $x\in S^1$ and
$t>0$, extends to a holomorphic motion of $S^1$ on the Riemann
sphere.

\vskip .2 cm

The direction (2)$\implies$(1) is elementary \cite[Proposition
2.1]{Sa}. The key idea for proving (1)$\implies$(2) is to show
that (1)$\implies$(3) by extending positive parameter $t$ into a
complex parameter with small imaginary part and showing that the
resulting map is a holomorphic motion gives (3)$\implies$(2). This
involved a geometric argument and the use of a result by
Ma\~{n}\'{e}, Sad and Sullivan \cite{MSS} on holomorphic motions.
Epstein, Marden and Markovic \cite{EMM} somewhat improved the
lower bound on the size of the imaginary part of the complex
parameter in their investigation of convex cores of plane domains.

\vskip .2 cm

At approximately the same time as in \cite{Sa}, Gardiner, Hu and
Lakic \cite{GHL} gave another proof of Theorem 1 using ordinary
differential equations. Their proof is more elementary than in
\cite{Sa} (since it does not use holomorphic motions) but it is
more computational and it does not show the analyticity of the
path $t\mapsto E_{t\mu}|_{S^1}$, as in \cite{Sa}, but only
differentiability of $t\mapsto E_{t\mu}|_{S^1}(x)$, $x\in S^1$.
Later on, another proof is given by Hu \cite{Hu} which is more
direct but used more computations as well. Both \cite{GHL},
\cite{Hu} consider stronger than quasisymmetric smoothness
properties of circle homeomorphisms.

\vskip .2cm

Since earthquakes are given by a geometric construction in a very
elementary terms on $\D$ (Thurston \cite[Introduction]{Th}), one
would hope to have a simple proof of (1)$\Longleftrightarrow$(2)
without much computations and using only elementary hyperbolic
geometry. We give a surprisingly easy proof of the equivalence of
the two conditions. The main difficulty is in proving
(1)$\implies$(2). The idea is to assume on the contrary that the
extension to $S^1$ of the earthquake map is not quasisymmetric. By
re-scalings of quadruples where the quasisymmetric constant tends
to infinity, we obtain a sequence of maps which cannot converge to
a homeomorphism. However, their corresponding earthquake measures
are the re-scalings of the bounded earthquake measure $\mu$. Hence
they are uniformly bounded and there exists a convergent
subsequence whose limit is a bounded earthquake measure. This
implies that the corresponding subsequence of earthquakes, when
properly normalized, converges to a homeomorphism of $S^1$. The
contradiction finishes the proof. The proof of the part
(2)$\implies$(1) from \cite[Theorem 1]{Sa} is briefly outlined for
the benefit of the reader.

\vskip .2 cm

An appropriate extension of the above method allows us to
characterize which earthquakes extend to symmetric maps of the
unit circle. We say that a bounded earthquake measure $\mu$ is
{\it asymptotically trivial} if $\sup_I\mu (I)\to 0$ as $\delta\to
0$, where the supremum is over all hyperbolic disks of a fixed
radius $r_0$ whose Euclidean distance to the boundary $S^1$ is at
most $\delta$.

\vskip .2 cm

\paragraph{\bf Theorem 2} {\it Let $h:S^1\to S^1$ be a
homeomorphisms and let $E_{\mu}:\D\to\D$ be the earthquake map
such that $E_{\mu}|_{S^1}=h$. Then the following are equivalent:}

(4) {\it $h$ is a symmetric map}

(5) {\it $\mu$ is asymptotically trivial earthquake measure. }

\section{Preliminaries}

\paragraph{\bf Definition 2.1} \cite{Th} Let $\lambda$ be a geodesic
lamination on $\D$. A {\it stratum} of $\lambda$ is either a
geodesic of $\lambda$ or a component of the complement of
$\lambda$ in $\D$, if any. An {\it earthquake} $E$ with the {\it
support} $\lambda$ is a surjective map $E:\D\to\D$ such that $E$
is a hyperbolic isometry when restricted to any stratum and, for
any two strata $A$ and $B$, the {\it comparison isometry}
$$
 E|_B\circ (E|_A)^{-1}
$$
is a hyperbolic translation whose axis weakly separates $A$ and
$B$, and which translates $B$ to the left as seen from $A$.

\vskip .2 cm

An earthquake $E$ of $\D$ continuously extends to a homeomorphism
of the boundary $\partial\D =S^1$. We denote by $E|_{S^1}$ the
extension. A remarkable fact is that any homeomorphism of $S^1$ is
obtained as a boundary extension of an earthquake \cite{Th}.

\vskip .2 cm

The translation length of $ E|_B\circ (E|_A)^{-1}$ is the first
approximation to the transverse measure deposited on a geodesic
arc $I$ connecting $A$ to $B$ corresponding to the earthquake $E$.
The {\it transverse measure} to $I$ is obtained by taking the
limit of the sum of translation lengths of comparison isometries
between finitely many consecutive strata intersecting $I$, as the
maximum distance between consecutive strata goes to zero. The
transverse measure to any geodesic arc is well-defined. Two
homeomorphisms of $S^1$ have the same transverse measures (to all
geodesic arcs) corresponding to their earthquakes if and only if
one homeomorphism is equal to the post-composition with a
hyperbolic isometry of the other homeomorphism (see \cite{Th}).
(Note that transverse measures determine earthquakes and circle
homeomorphisms uniquely up to a post-composition by a hyperbolic
isometry of $D$. Thus, given a transverse measure, we need to
normalize the earthquake by fixing the hyperbolic isometry.) The
transverse measures to arcs are invariant under isotopies of $\D$
which preserve the support geodesic lamination $\lambda$. The
family of the transverse measures to all geodesic arcs in $\D$ is
called an {\it earthquake measure}.

\vskip .2 cm

Let $\mathcal{G}$ be the set of all unoriented geodesics of $\D$.
Then $\mathcal{G}$ is homeomorphic to $(S^1\times S^1-diag)/\sim$,
where $(a,b)\sim(b,a)$. It is sometimes useful to think of
transverse measures to the support geodesic laminations $\lambda$
as measures on the space of geodesics $\mathcal{G}$ whose support
is a geodesic lamination $\lambda$. An {\it earthquake measure}
$\mu$ is a positive Radon measure on $\mathcal{G}$ whose support
is a geodesic lamination.

\vskip .2 cm

Let $I$ be an arbitrary geodesic arc of length $1$. Then $\mu (I)$
is the total measure of $I$ deposited by an earthquake measure
$\mu$, or alternatively, $\mu (I)$  is the $\mu$-mass of all
geodesics in $\D$ intersecting $I$. Thurston \cite{Th} introduced
the {\it norm} of an earthquake measure $\mu$ by
$$
\|\mu\| :=\sup_{I}\mu (I),
$$
where the supremum is over all closed geodesic arcs $I$ of length
$1$. An earthquake measure $\mu$ is said to be {\it bounded} if
$\|\mu \|<\infty$.

\vskip .2 cm

Let $r_0>0$ be a fixed number smaller than the hyperbolic radius
of inscribed disk in an ideal hyperbolic triangle in $\D$. If $D$
is a hyperbolic disk with radius $r_0$, then the set of geodesic,
from an arbitrary geodesic lamination of $\D$, which intersect $D$
also intersect a single hyperbolic arc of length $2r_0$. We define
the measure of $D$ deposited by an earthquake measure $\mu$ of
$\D$ to be equal to the measure deposited to the arc which is
transverse to the subset of support of $\mu$ which intersects $D$.

\vskip .2 cm

Let $D$ be an arbitrary hyperbolic disk in $\D$. Define $\delta
(D)$ to be the Euclidean distance of $D$ to the boundary $S^1$. A
bounded earthquake measure $\mu$ is {\it asymptotically trivial}
if
$$
\sup_{D,\ \delta (D)\leq t}\mu (D)\to 0
$$
as $t\to 0$, where the supremum is over all hyperbolic disks $D$
with radius $r_0$ whose distance to the boundary $\delta (D)$ is
at most $t$.

\vskip .2 cm

An {\it earthquake cocycle} with support geodesic lamination
$\lambda$ is a map $E:\D\times\D \to PSL_2(\mathbb{R})$ with the
cocycle property such that, given two strata $A$ and $B$ for
$\lambda$, we have $E(z_1,z_2)=E(w_1,w_2)$ whenever $z_1,w_1\in A$
and $z_2,w_2\in B$. In addition, $E(z_1,z_2)$ is required to be a
hyperbolic translation with the axis weakly separating the strata
$A$ and $B$ containing $z_1$ and $z_2$ and translating $B$ to the
left as seen from $A$. Given an earthquake measure $\mu$ there
exists a corresponding earthquake cocycle which defines a
piecewise isometric, injective map $E:\D\to\D$ \cite{EpM}. (If $E$
is surjective, then it is an earthquake map.) The earthquake
obtained in this way has its measure equal to $\mu$ \cite{EpM},
\cite{GHL}, \cite{Sa}.

\section{The weak convergence}

Denote by $\mathcal{G}_z$, for $z\in\D$, the set of geodesics in
$\D$ which contain $z$. Given $z,w\in\D$ denote by $[z,w]$ the
geodesic arc in $\D$ between $z$ and $w$. If $K$ is a subset of
$\D$, denote by $\mathcal{G}_K$ the set of geodesics of $\D$ which
intersect $K$. The following lemma is essentially proved by
Epstein and Marden \cite[Theorem 3.11.5]{EpM}. Their statement is
for sequences of finite complex earthquake measures, but it
immediately extends to arbitrary (not necessarily finite)
sequences of positive earthquakes with a simpler proof.

\vskip .2 cm

\paragraph{\bf Lemma 3.1} \cite{EpM} {\it Let $\mu_i,\mu$ be earthquake
measures such that $\mu_i\to\mu$ as $i\to\infty$ in the weak*
topology. Then the sequence $E_{\mu_i}(z_1,z_2)$, for any
$z_1,z_2\in\D$, of earthquake cocycles has a convergent
subsequence. Any limit $E(z_1,z_2)$ of a subsequence
$E_{\mu_{i_j}}(z_1,z_2)$ satisfies
$$
E_{\mu}(z_1,z_2)=T_{g_2}^{a_2}\circ E(z_1,z_2)\circ T_{g_1}^{a_1},
$$
where $a_k=\mu(g_k)-\nu(g_k)$, $\nu$ is the weak* limit of
$\mu_{i_j}$ restricted to $\mathcal{G}_{[z_1,z_2]}$, and
$T_{g_k}^{a_k}$, $k=1,2$, is either the identity if $a_k=0$, or it
is a hyperbolic translation with the axis $g_k$ ($z_k\in g_k$) in
the support of $\mu$ and the translation length $a_k$. }

\vskip .2 cm

The following lemma is to be expected, but the proof is somewhat
subtle. Note that we require the sequence of earthquake measures
$\mu_i$ to be uniformly bounded.

\vskip .2 cm

\paragraph{\bf Lemma 3.2} {\it Let $\mu_i,\mu$, for $i=1,2,\ldots$,
be uniformly bounded earthquake measures (i.e. $\|\mu\| ,\|\mu_i\|
\leq M<\infty$ for all $i$) on the unit disk $\D$. If
$\mu_i\to\mu$ as $i\to\infty$ in the weak* topology then, for each
$x\in S^1$,
$$
E_{\mu_i}|_{S^1}(x)\to E_{\mu}|_{S^1}(x)
$$
as $i\to\infty$, when the earthquakes
$E_{\mu_i}|_{S^1},E_{\mu}|_{S^1}$ are properly normalized.}

\vskip .2 cm

\paragraph{\bf Proof} Recall that $E_{\mu_i}|_{S^1},E_{\mu}|_{S^1}$ are
well-defined up to post-composition by a hyperbolic isometry of
$\D$. We normalize the earthquake maps as follows. Either the
earthquake measure $\mu$ has a complementary gap $A$ in $\D$ or
the support of $\mu$ foliates $\D$. In the first case, we fix
$E_{\mu}$ to be the identity on the gap $A$. In the second case,
we fix a geodesic $l$ in the support of $\mu$ such that $\mu(l)=0$
and set $E_{\mu}|_l=id$.

\vskip .2 cm

If the support of $\mu$ has a gap $A$ then either the support of
$\mu_i$ has a gap $A_i$ which intersects $A$ with at least one
boundary side contained in a fixed compact subset of $\mathcal{G}$
or there exists a geodesic $l_i$ in the support of $\mu_i$ of zero
$\mu_i$-measure which intersects $A$ and is contained in the fixed
compact subset of $\mathcal{G}$ (for $i$ large enough depending on
the fixed compact subset of $\mathcal{G}$). If the support of
$\mu$ has no gaps then there exists a geodesic $l_i$ in the
support of $\mu_i$ such that $\mu_i(l_i)\to 0$ and $l_i\to l$ as
$i\to\infty$. Then we set either $E_{\mu_i}|_{A_i}=id$ or
$E_{\mu_i}|_{l_i}=id$.

\vskip .2 cm

We fix $r$, $0<r<1$, and define $\mu'$ to be the restriction of
$\mu$ to the set of geodesics which intersect the euclidean disk
$\D_r$ centered at $0$ of radius $r$. We think of $\mu'$ as a new
earthquake measure on $\G$ whose support geodesics intersect
$\D_r$. Then there exists a sequence $r_i\to r$, $r_i\geq r$, such
that $\mu_i'\to\mu'$ as $i\to\infty$, where $\mu_i'$ is the
restriction of $\mu_i$ to the subset of $\mathcal{G}$ whose
geodesics intersect $\D_{r_i}$. For $r$ large enough, either $A$
is contained in or is equal to a stratum $A'$ of $\mu'$, or $l$ is
in the support of $\mu'$. Then, for $i$ large enough, either $A_i$
is contained in a stratum $A_i'$ of $\mu_i'$, or $l_i$ is in the
support of $\mu_i'$. We normalize $E_{\mu'}$ and $E_{\mu_i'}$ in a
corresponding manner.

\vskip .2 cm

We show that $E_{\mu_i'}|_{S^1}(x)\to E_{\mu'}|_{S^1}(x)$ as
$i\to\infty$, for all $x\in S^1$. Since the supports of $\mu'$ and
$\mu_i'$ are contained in a compact subset of $\mathcal{G}$, each
$x\in S^1$ is on the boundary of at least one strata of the
measures $\mu'$ and $\mu_i'$. Note that $E_{\mu'}$ has at most
countably many leaves in the support with non-zero measure which
intersect $S^1$ in at most countably many points. It is enough to
prove the convergence outside these points of $S^1$ due to the
fact that $E_{\mu'}|_{S^1}$ and $E_{\mu_i'}|_{S^1}$ are order
preserving maps. If $x\in S^1$ is a point which lies on the
boundary of a stratum $B$ of $\mu'$ (if the stratum is a geodesic
then the $\mu$-measure is zero by our assumption) then the
measures $\mu_i'$ restricted to the geodesic arc $I$ connecting
the fixed stratum ($A'$ or $l'$) of $\mu'$ to the stratum $B$ ($I$
connects the interior of the strata $A'$ with the interior of the
strata $B$ if the strata are gaps) converge to the restriction of
the measure $\mu'$ to $I$. Thus the earthquake cocycles
$E_{\mu_i'}$ corresponding to the endpoints of $I$ converge to
$E_{\mu'}$ by Lemma 3.1. Moreover, since each geodesic in the
support of $\mu'$ is approximated by the geodesics of the support
of $\mu_i'$, there exists a sequence $I_i$ of closed geodesic arcs
(which lie on the same ideal geodesic as $I$) with one endpoint at
$A_i'$ or $l_i'$ and the other endpoint at $B_i$ (where $B_i$ is a
stratum of $\mu_i'$ converging to $B$ and containing $x\in S^1$ on
its boundary) which converge toward $I$. It follows by our choice
of $B_i$ that the restriction of $\mu_i'$ to $I_i$ converges to
the restriction of $\mu'$ to $I$ (when both restrictions are
considered as measures on $\mathcal{G}$) in the weak* topology. By
Lemma 3.1, we get that the cocycles for $\mu_i'$ with respect to
the endpoints of $I_i$ converge to the cocycle for $\mu'$ with
respect to $I$ when they are considered as hyperbolic isometries
of $\D$. This implies the desired convergence for any fixed $r<1$.

\vskip .2 cm

Therefore, to show that $E_{\mu_i}|_{S^1}(x)\to E_{\mu}|_{S^1}(x)$
as $i\to\infty$, it is enough to show that
$E_{\mu_i'}|_{S^1}(x)\to E_{\mu_i}|_{S^1}(x)$ and
$E_{\mu'}|_{S^1}(x)\to E_{\mu}|_{S^1}(x)$ as $r\to 1$ independent
of $i$. By the cocycle property, we have
$E_{\mu_i}|_{S^1}(x)=E_{\mu_i'}|_{S^1}\circ
E_{\mu_i''}|_{S^1}(x)$, where $E_{\mu_i''}$ is normalized to be
the identity on the stratum $C$ which contains $\D_r$ and
$\mu_i''=\mu_i -\mu_i'$. However
$E_{\mu_i}|_{S^1}(x)=E_{\tilde{\mu}_i''}|_{S^1}\circ
E_{\mu_i'}|_{S^1}(x)$, where
$\tilde{\mu}_i''=(E_{\mu_i'}|_{S^1})^{*}(\mu_i'')$ and
$E_{\tilde{\mu_i}''}$ is normalized to be the identity on the
stratum $E_{\mu_i'}(C)$. Note that the support of $\mu_i''$ does
not intersect $\D_r$ except possibly on the boundary
$\partial\D_r$ which implies that the diameter of each stratum of
$\mu_i''$ outside $\D_r$ goes to zero as $r\to 1$. The image of
the support of $\mu_i''$ under $E_{\mu_i'}$ is not increasing by
much, namely its diameter also goes to zero as $r\to 1$. This
follows by an easy observation that bounded earthquakes decrease
the distance between their support geodesics by a bounded amount
\cite{Th}, \cite[Lemma 2.1, Proposition 2.1]{Sa}. Since the
distance between $E_{\mu_i'}|_{S^1}(x)$ and $E_{\mu_i}|_{S^1}(x)$
is at most the diameter of the support of $\tilde{\mu}_i''$, we
get the desired convergence. $\Box$

\vskip .2 cm

In the following proposition we show that the converse is true.
This is the key point in our proof of Theorem 1.

\vskip .2 cm

\paragraph{\bf Proposition 3.3} {\it Let $\mu ,\mu_i$ be uniformly bounded
earthquake measures on $\D$. Then $\mu_i\to\mu$ in the weak*
topology as $i\to\infty$ if and only if there exist normalization
of earthquake maps $E_{\mu_i}|_{S^1},E_{\mu}|_{S^1}$ such that
$E_{\mu_i}|_{S^1}(x)\to E_{\mu}|_{S^1}(x)$ for each $x\in S^1$, as
$i\to\infty$.}

\vskip .2 cm

\paragraph{\bf Proof} The ``if'' part is proved in Lemma 3.2 above.
We show the ``only if'' part.

\vskip .2 cm

We suppose that $E_{\mu_i}|_{S^1}\to E_{\mu}|_{S^1}$ pointwise as
$i\to\infty$. Assume on the contrary that $\mu_i$ does not
converge to $\mu$ in the weak* topology. By the compactness of the
space of probability measures on a compact space, there exists a
subsequence $\mu_{i_j}$ of $\mu_i$ which converges to a measure
$\nu$ in the weak* topology (by the standard diagonal argument).
It is clear that $\nu$ is an earthquake measure bounded by the
same constant as the sequence $\mu_i$. Our assumption implies that
$\nu\neq\mu$. By Lemma 3.2, there exist normalization of
$E_{\mu_{i_j}}$ and $E_{\nu}$ such that $E_{\mu_{i_j}}|_{S^1}\to
E_{\nu}|_{S^1}$ pointwise, as $j\to\infty$. However, since
$E_{\mu_i}|_{S^1}\to E_{\mu}|_{S^1}$ with possibly different
normalization and since two normalizations differ by a
post-composition with a hyperbolic isometry of $\D$, we conclude
that $E_{\mu}|_{S^1}=\gamma\circ E_{\nu}|_{S^1}$ for some
hyperbolic isometry $\gamma$. By the uniqueness of the earthquake
measure \cite{Th}, we have that $\nu =\mu$ which gives a
contradiction. Therefore $\mu_i\to\mu$. $\Box$

\section{Proofs of theorems}

\paragraph{\bf Proof of Theorem 1} We sketch the proof of
(2)$\implies$(1) from \cite{Sa}. Assume on the contrary that $\mu$
is unbounded. If there is a sequence of geodesic $l_i$ in the
support of $\mu$ whose measure $\mu (l_i)\to\infty$ then we find a
sequence of quadruples with fixed cross-ratios whose images have
unbounded cross-ratios. It is enough to take quadruples which have
one point in each half-plane of $\D -l_i$ and the other two points
to be the endpoints of $l_i$. Since the earthquake translates to
the left and $\mu (l_i)\to\infty$, it follows that the
cross-ratios of the images are unbounded. Thus the boundary map is
not quasisymmetric. This is a contradiction. In the case that
there is no such $l_i$ in the support of $\mu$, then we can find a
sequence of subset $\mathcal{L}_i$ of the support of $\mu$ whose
each geodesic has one endpoint in an interval $(a_i,b_i)$ and the
other endpoint in an interval $(c_i,d_i)$ with the cross-ratio of
$(a_i,b_i,c_i,d_i)$ converging to infinity. Therefore,
$\mathcal{L}_i$ is close to being a single geodesic and similar
argument applies.

\vskip .2 cm

We prove (1)$\implies$(2). Assume on the contrary that
$h:=E_{\mu}|_{S^1}$ is not quasisymmetric. We already know that
$h$ is a homeomorphism by \cite{Th}, \cite[Proposition 2.1]{Sa},
or \cite{GHL}. Then there exists a sequence $(a_i,b_i,c_i,d_i)$ of
counterclockwise oriented quadruplets of points on $S^1$ with
cross-ratio $2$ such that $cr(h(a_i),h(b_i),h(c_i),h(d_i))\to
\infty$ as $i\to\infty$. Let $A_i$ be the hyperbolic isometry such
that $A_i:(a_i,b_i,c_i,d_i)\mapsto (1,i,-1,-i)$. Then $h_i:=h\circ
A_i^{-1}$ is a sequence of maps which are not quasisymmetric. Note
that $\gamma_i\circ h_i=E_{\mu_i}|_{S^1}$ for some hyperbolic
isometry $\gamma_i$, where $\mu_i:=(A_i)^{*}(\mu )$. Since $\mu_i$
is a sequence of uniformly bounded earthquake measures, there is a
subsequence $\mu_{i_j}$ which converges to a bounded earthquake
measure $\sigma$. By Proposition 3.3, we get that
$E_{\mu_{i_j}}|_{S^1}\to E_{\sigma}|_{S^1}$ pointwise as
$j\to\infty$ when the earthquakes are properly normalized.

\vskip .2 cm

On one hand, $E_{\mu_{i_j}}|_{S^1}=\gamma_{i_j}\circ h_{i_j}$
where $\gamma_{i_j}$ is a hyperbolic isometry of $\D$. On the
other hand, $E_{\mu_{i_j}}|_{S^1}$ converges to a homeomorphism of
$S^1$ which implies that the cross-ratios of
$E_{\mu_{i_j}}|_{S^1}(1,i,-1,-i)$ are uniformly bounded. But this
is in the contradiction with $cr(h(a_i),h(b_i),h(c_i),h(d_i))\to
\infty$ as $i\to\infty$ because
$cr(h(a_{i_j}),h(b_{i_j}),h(c_{i_j}),h(d_{i_j}))=
cr(h_{i_j}(1),h_{i_j}(i),h_{i_j}(-1),h_{i_j}(-i)).$ Thus $h$ is
quasisymmetric. $\Box$

\vskip .4 cm

\paragraph{\bf Proof of Theorem 2} We first show (4)$\implies$(5) part.
Assume on the contrary that $\mu$ is not asymptotically trivial
earthquake measure. Then there exists a sequence of hyperbolic
disks $D_i$ whose radii are $r_0$ and whose centers are
$d_i\in\D$, $|d_i|\to 1$, such that $\mu (D_i)\geq m>0$. Let $A_i$
be the hyperbolic translation with the axis the radius of $\D$
through $d_i$ which maps $d_i$ onto $0$. Let
$$
\mu_i:=A_i^{*}(\mu )
$$
and let $D$ be the hyperbolic disk with the center $0$ and the
radius $r_0$.

\vskip .2 cm

Since $\|\mu_i\| =\|\mu\|$, for each $i$, it follows that the
sequence of positive earthquake measures $\mu_i$ is uniformly
bounded. Then there exists a convergent subsequence
$\mu_{i_j}\to\sigma$ as $j\to\infty$ in the weak* topology, where
$\sigma$ is a bounded earthquake measure. Note that $\sigma
(D)\geq m>0$ because $\mu_{i_j}(D)=\mu (D_{i_j})\geq m$. By
Proposition 3.3, there exists a normalization of earthquakes
$E_{\sigma},E_{\mu_{i_j}}$ such that $E_{\mu_{i_j}}|_{S^1}\to
E_{\sigma}|_{S^1}$ pointwise as $j\to\infty$. In fact, it follows
from the proof that we can normalize $E_{\sigma}$ to be the
identity on a stratum intersecting $D$. Since $\sigma\neq 0$, we
have that $h^{*}:=E_{\sigma}|_{S^1}$ is not the restriction to
$S^1$ of a hyperbolic isometry of $\D$ by the uniqueness of the
earthquake measures \cite{Th}. Define
$h_{i_j}:=E_{\mu_{i_j}}|_{S^1}$.

\vskip .2 cm

Let $f^{*}:=ex(h^{*})$, $f_{i_j}:=ex(h_{i_j})$ and $f:=ex(h)$ be
barycentric extensions (see \cite{DE}) of $h^{*}$, $h_{i_j}$ and
$h$, respectively. Since $h$ is symmetric, it follows that $f$ is
asymptotically conformal, i.e. given any $\epsilon >0$ there
exists a compact subset $K$ of $\D$ such that the supremum of the
absolute value of the Beltrami coefficient of $f$ over $\D -K$ is
at most $\epsilon$. Note that $f_{i_j}\to f^{*}$ pointwise by the
properties of barycentric extension \cite{DE} because $h_{i_j}\to
h^{*}$. Moreover, since $h_{i_j}$ is obtained as an extension to
the boundary of the earthquake $E_{\mu_{i_j}}$ whose measure is
the push-forward of the measure for $h$ by the hyperbolic isometry
$A_{i_j}$, we get that
$$h_{i_j}=B_{i_j}\circ h\circ A_{i_j}^{-1},$$
where $B_{i_j}$ is a hyperbolic isometry. Since barycentric
extension is conformaly natural, we get
$$
f_{i_j}=B_{i_j}\circ f\circ A_{i_j}^{-1}.
$$
By taking Beltrami coefficients of the left and of the right side
in the above equation, we get that
$$
Belt(f_{i_j})=Belt(f)\circ
A_{i_j}^{-1}\frac{\overline{A_{i_j}^{-1}}}{A_{i_j}^{-1}}.
$$
Since $f$ is asymptotically conformal and by the definition of
$A_{i_j}$, we conclude that $Belt(f)\circ A_{i_j}^{-1}\to 0$ as
$j\to\infty$ uniformly on compact subset of $\D$. Thus $f^{*}$ is
a hyperbolic isometry of $\D$ and its extension $h^{*}$ to the
boundary is necessarily a hyperbolic isometry. Contradiction.
Therefore $\mu$ is asymptotically trivial.

\vskip .3 cm

We show (5)$\implies$(4) part. Assume on the contrary that
$h=E_{\mu}|_{S^1}$ is not symmetric. We recall a characterization
of symmetric maps from \cite{Su}. Cover $S^1$ by finitely many
charts which are closed, bounded intervals in $\mathbb{R}$. In
each chart, a 4-tuple $(a,b,c,d)$, $a,b,c,d\in \mathbb{R}$, is
called {\it standard} if $|a-b|=|b-c|=|c-d|$. A homeomorphism of
$S^1$ is symmetric if and only if, for any standard 4-tuple
$(a,b,c,d)$ in any chart, the {\it cross-ratio distortion}
$\log\frac{cr(h(a,b,c,d))}{cr(a,b,c,d)}$ of $h$ converges to $0$
as $|a-b|\to 0$. (Any standard 4-tuple on $\mathbb{R}$ has
cross-ratio $4/3$. We choose the chart maps for $S^1$ to be
M\"obius. Then the standard 4-tuples in charts are the images of
4-tuples on $S^1$ with fixed cross-ratios $4/3$. By abuse of
notation, we identify standard 4-tuples in charts with their
images on $S^1$ under chart maps.)

\vskip .2 cm

Since $h$ is not symmetric, there exists a sequence of standard
4-tuples $(a_i,b_i,c_i,d_i)$, with $|a_i-b_i|\to 0$, in a fixed
chart such that $\log
\frac{cr(h(a_i,b_i,c_i,d_i))}{cr(a_i,b_i,c_i,d_i)}$ does not
converge to $0$. Let $A_i$ be the hyperbolic isometry of $\D$
which maps $(a_i,b_i,c_i,d_i)$ onto $(-1,a,b,1)$ with
$cr(-1,a,b,1)=4/3$. Define $\mu_i:=A_i^{*}(\mu )$. Note that
$\mu_i\to 0$ as $i\to\infty$ in the weak* topology because $\mu$
is asymptotically trivial. Then $E_{\mu_i}|_{S^1}\to id$ pointwise
by Proposition 3.3, when $E_{\mu_i}|_{S^1}$ are properly
normalized. Define $h_i:=h\circ A_i^{-1}$. Then
$E_{\mu_i}|_{S^1}=B_i\circ h_i$, for some hyperbolic isometry
$B_i$. Since $B_i\circ h_i\to id$ pointwise, we conclude that
$\log \frac{cr(B_i\circ h_i(-1,a,b,1))}{cr(-1,a,b,1)}\to 0$ as
$i\to\infty$. On the other hand, $\log \frac{cr(B_i\circ
h_i(-1,a,b,1))}{cr(-1,a,b,1)}=\log
\frac{cr(h_i(-1,a,b,1))}{cr(-1,a,b,1)}=\log \frac{cr(
h(A_i^{-1}(-1,a,b,1))}{cr(-1,a,b,1)}=\log \frac{cr(
h(a_i,b_i,c_i,d_i))}{cr(-1,a,b,1)}$ does not converge to $0$ by
the above. Contradiction. Therefore $h$ is symmetric. $\Box$

\end{document}